\documentclass[11pt,twoside]{amsart}

\usepackage{bbm,amsmath,amssymb,amsthm,amscd,latexsym,mathrsfs}
\usepackage[all]{xy}
\usepackage{pb-diagram,pb-xy}
\usepackage[pdftex]{hyperref}

\newtheorem{Prop}{Proposition}[subsection]
\newtheorem{Thm}[Prop]{Theorem}
\newtheorem{Theorem}{Theorem}

\newtheorem{Lemma}[Prop]{Lemma}
\newtheorem{Cor}[Prop]{Corollary}

\newenvironment{Proof}{\par \noindent \textit{Proof.}}{\hfill $\Box$ \par}

\theoremstyle{definition}
\newtheorem{Def}[Prop]{Definition}

\newcommand{\OO}{\mathcal{O}}

\newcommand{\CC}{\mathbb{C}}
\newcommand{\ZZ}{\mathbb{Z}}
\newcommand{\HH}{\mathbb{H}}
\newcommand{\RR}{\mathbb{R}}
\newcommand{\NN}{\mathbb{N}}
\newcommand{\QQ}{\mathbb{Q}}
\newcommand{\PP}{\mathbb{P}}

\newcommand{\UU}{\mathrm{U}}

\DeclareMathOperator{\Hom}{Hom}

\DeclareMathOperator{\im}{im}

\DeclareMathOperator{\NS}{NS}

\DeclareMathOperator{\Res}{Res}
\DeclareMathOperator{\Hdg}{Hdg}
\DeclareMathOperator{\End}{End}
\DeclareMathOperator{\id}{id}
\DeclareMathOperator{\cores}{Cores}
\DeclareMathOperator{\Spin}{Spin}

\DeclareMathOperator{\SO}{SO}
\DeclareMathOperator{\SMT}{SMT}

\DeclareMathOperator{\cc}{c}

\begin{document}

\title{The Hodge conjecture for self-products of certain K3 surfaces}
\author{Ulrich Schlickewei}

\address{Mathematisches Institut der Universit{\"a}t Bonn, Endenicher Allee 60,
53115 Bonn, Germany} 
\email{uli@math.uni-bonn.de}

\begin{abstract}
We use a result of van Geemen \cite{vG2} to determine the endomorphism algebra of the Kuga--Satake variety
of a K3 surface with real multiplication. This is applied to prove the Hodge conjecture for self-products of double covers
of $\PP^2$ which are ramified along six lines.
\end{abstract}

\maketitle
\let\thefootnote\relax\footnotetext{This work was supported by the SFB/TR 45 `Periods,
Moduli Spaces and Arithmetic of Algebraic Varieties' of the DFG
(German Research Foundation) and by the Bonn International Graduate School in Mathematics (BIGS).}

\begin{section}{Introduction}
Let $S$ be a complex K3 surface, i.e.\ a smooth, projective surface over $\CC$ satisfying
$H^1(S,\OO_S) = 0$ and $\omega_S \simeq \OO_S$. Let $T(S) \subset H^2(S,\QQ)$ be the
rational transcendental lattice of $S$, defined as the orthogonal complement of
the N\'eron--Severi group with respect to the intersection form.
The algebra $E_S := \End_{\Hdg}(T(S))$ of endomorphisms of $T(S)$ which preserve the Hodge decomposition can
be interpreted as a subspace of the space of (2,2)-classes on the self-product $S \times S$.
The Hodge conjecture for $S \times S$ predicts that $E_S$ consists of linear combinations
of fundamental classes of algebraic surfaces in $S \times S$. Using the Lefschetz theorem on (1,1)-classes,
it is easily seen that conversely the Hodge conjecture for $S \times S$ holds if $E_S$ is
generated by algebraic classes.

\vspace{2ex}
Mukai \cite{Mu} used his theory of moduli spaces
of sheaves to prove that if the Picard number of $S$ is at least
11, then any $\varphi \in E_S$ which preserves the intersection
form on $T(S)$ can be represented as a linear combination of fundamental classes of algebraic cycles. Later
this result was improved by Nikulin \cite{N} on the base of lattice-theoretic arguments to the case that the 
Picard number of $S$ is at least 5. In \cite{Mu2}, Mukai announced that using the theory of moduli spaces 
of twisted sheaves, the hypothesis on the Picard number could be omitted.

But how many isometries do exist in the algebra $E_S$?
Results of Zarhin \cite{Zarhin} imply that $E_S$ is an algebraic number field,
which is either totally real (we say that $S$ has \emph{real multiplication})
or a CM field ($S$ has \emph{complex multiplication}). Isometries of $T(S)$ correspond to 
elements of norm 1 in $E_S$.
If $S$ has complex multiplication, one can use the 
fact that CM fields are generated as $\QQ$-vector spaces by elements of norm 1 to see that $E_S$ is generated by 
isometries. In combination with Mukai's results, this proves the Hodge conjecture for
self-products of K3 surfaces with complex multiplication and with Picard number at least 5. This 
was noticed by Ram\'on-Mar\'i \cite{Ma}. If $S$ has real multiplication, the only Hodge isometries in $E_S$ are plus or minus the
identity. Thus, Mukai's results are no longer sufficient to prove the
algebraicity of interesting classes in $E_S$.

\vspace{2ex}
In order to approach the case of real multiplication one passes from
K3 surfaces to Abelian varieties by associating to a K3 surface
$S$ its Kuga--Satake Abelian variety $A$. By construction, see \cite{KS},
there exists an inclusion of Hodge structures $T(S)\subset H^2(A\times A,\QQ)$.
Van Geemen \cite{vG2}
studied the Kuga--Satake variety of a K3 surface
with real multiplication. He discovered that the corestriction of a certain Clifford algebra
over $E_S$ plays an important role for the Kuga--Satake variety of $S$. We rephrase and slightly improve
his result which then reads as follows:

\begin{Theorem} \label{ThmZerlegungKS}
 Let $S$ be a K3 surface with real multiplication by a totally real number field $E_S$
of degree $d$ over $\QQ$. Let $A$ be a Kuga--Satake variety of $S$. 

Then there exists an Abelian variety $B$ such that $A$ is isogenous to $B^{2^{d-1}}$. The endomorphism algebra of $B$ is
\begin{equation*}
  \End_{\QQ}(B) = \cores_{E /\QQ} C^0(Q).
\end{equation*}
\end{Theorem}

Here, $Q: T \times T \to E_S$ is a quadratic form on $T$ which already appeared 
in Zarhin's paper \cite{Zarhin} and which will be reintroduced in Section \ref{SectionSMTZarhin},
$C^0(Q)$ is the even Clifford algebra of $Q$ over $E_S$ and $\cores_{E/\QQ}C^0(Q)$ is
the corestriction of this algebra. 
The corestriction of algebras will be reviewed in Section \ref{SectionCorestriction}.

\vspace{2ex}
Theorem \ref{ThmZerlegungKS} leads to a better understanding of the phenomenon of real multiplication
for K3 surfaces 
by allowing us to calculate the endomorphism algebra of the corresponding Kuga--Satake varieties. 
However, since the Kuga--Satake construction is purely
Hodge-theoretic, this still gives no geometric explanation. Therefore, we focus on one of the 
few families of K3 surfaces for which the Kuga--Satake correspondence
has been understood geometrically. This is the family of double covers of $\PP^2$ ramified along six lines.
Paranjape \cite{P} found an explicit cycle on $S \times A \times A$ which realizes the inclusion of Hodge
structures $T(S) \subset H^2(A \times A,\QQ)$. Building on the decomposition result for Kuga--Satake varieties we derive

\begin{Theorem} \label{ThmHCDC}
 Let $S$ be a K3 surface which is a double cover of $\PP^2$ ramified along six lines. Then the Hodge
 conjecture is true for $S \times S$.
\end{Theorem}

As pointed out by van Geemen \cite{vG2}, there are one-dimensional sub-families 
of the family of such double covers with real multiplication
by a totally real quadratic number field. In conjunction with our
Theorem \ref{ThmHCDC}, this allows us to produce examples of K3 surfaces $S$
with non-trivial real multiplication for which $\End_{\Hdg}(T(S))$ is generated by algebraic classes.
We could not find examples of this type in the existing literature.

\vspace{2ex}
The \emph{plan of the paper} is as follows: In Section \ref{SectionHS} we review Zarhin's results on the endomorphism algebra 
and on the special Mumford--Tate group of an
irreducible Hodge structure of K3 type. Also, we recall from \cite{vG2} how a Hodge structure of K3 type with real multiplication 
splits over a finite extension of $\QQ$.

Section \ref{SectionKS} is devoted to the proof of Theorem \ref{ThmZerlegungKS}. After reviewing the definition
of the corestriction of algebras, we explain in detail how the Galois group of a normal closure of $E_S$ acts on the
Kuga--Satake Hodge structure. This is the key of the proof.

In the final Section \ref{SectionDoubleCovers} we study double covers of $\PP^2$ ramified along six lines. We 
review results of Lombardo \cite{L} on the Kuga--Satake variety of such K3 surfaces, of Schoen \cite{S} and van Geemen
\cite{vG3} on the Hodge conjecture for certain Abelian varieties of Weil type and of course Paranjape's \cite{P} result on the
algebraicity of the Kuga--Satake correspondence. Together with Theorem \ref{ThmZerlegungKS}, they
lead to the proof of Theorem \ref{ThmHCDC}.

\vspace{2ex} \noindent
\emph{Acknowledgements.} This work is a part of my Ph.D.\ thesis prepared at the University of Bonn. It is a
great pleasure to thank my advisor Daniel Huybrechts for suggesting this interesting topic and for constantly
supporting me.

During a four week stay at the University of Milano I had many fruitful discussions with Bert
van Geemen. I am most grateful to him for his insights.
\end{section}

\begin{section}{Hodge structures of K3 type with 
real multiplication} \label{SectionHS}
\begin{subsection}{Hodge structures of K3 type and their endomorphisms} 
 Let $\mathrm{U}(1)$ be the one-dimensional unitary group which is a real algebraic group. 
 To fix notations we recall that a Hodge structure of weight $k$ is a finite-dimensional $\QQ$-vector space $T$ together
 with a morphism of real algebraic groups $h: \mathrm{U}(1) \to \mathrm{GL}(T)_{\RR}$ such that 
 for $z \in \mathrm{U}(1)(\RR)\subset \CC$ 
 the $\CC$-linear extension of the 
 endomorphism $h(z)$ is diagonalizable with eigenvalues $z^p \overline{z}^q$ where $p+q= k$
 and $p,q \ge 0$ 
 (cf.\ e.g.\ 
 \cite[1.1]{vG2}). The eigenspace to $z^p \overline{z}^q$ is denoted by $T^{p,q} \subset T_{\CC}$.
 
 A polarization of a weight $k$ Hodge structure $(T,h)$ is a bilinear form $q: T \times T \to \QQ$ which is 
 $\mathrm{U}(1)$-invariant and which has the property that 
 $(-1)^{k (k-1)/2}q(*,h(i)*): T_{\RR} \times T_{\RR} \to T_{\RR}$ is a symmetric,
 positive definite bilinear form.

\begin{Def} A \emph{Hodge structure of K3 type $(T,h,q)$} consists of a $\QQ$-Hodge structure $(T, h: \mathrm{U}(1) \to \mathrm{GL}(T)_{\RR})$ of weight 2 with 
 $\dim_{\CC} T^{2,0} = 1$ together
 with a polarization $q : T \times T \to \QQ$. 
\end{Def}

 \noindent \emph{Examples.} The second primitive (rational) cohomology and the (rational) transcendental lattice 
 of a projective K3 surface yield examples of Hodge structures of K3
 type. More generally, the second primitive cohomology and the Beauville--Bogomolov orthogonal complement of the 
 N\'eron--Severi group
 of an irreducible symplectic variety are Hodge structures of K3 type \cite[Part III]{GHJ}.
 
 \vspace{2ex}

 Consider the Hodge decomposition 
 \begin{equation*}
   T_{\CC} := T \otimes_{\QQ} \CC = T^{2,0} \oplus T^{1,1} \oplus T^{0,2}.
 \end{equation*}
 Since the quadratic form $q$ is a polarization, this decomposition is $q$-orthogonal. 
 Moreover, 
 $q$ is positive definite on $(T^{2,0} \oplus T^{0,2})\cap T_{\RR}$ and negative definite on $T^{1,1} \cap T_{\RR}$.

 \vspace{2ex}
 Assume that $T$ is an irreducible Hodge structure. Let $E:= \End_{\Hdg}(T)$ be the division algebra of endomorphisms of Hodge 
 structures of $T$. Let $':E \to E$ be the involution given by adjunction with respect to $q$ and let $E_0 \subset E$ be the
 subalgebra of $E$ formed by $q$-self-adjoint endomorphisms.
 \begin{Thm}[Zarhin \cite{Zarhin}] \label{ThmZarhin1} The map 
 \begin{equation*} 
      \epsilon: E \to \CC, \; \; \; e \mapsto \text{eigenvalue of} \; e \; \text{on the eigenspace} \; T^{2,0}
 \end{equation*}
 identifies $E$ with a subfield of $\CC$. Moreover, $E_0$ is a totally real number field and the following two cases
 are possible:

 $\bullet$ $E_0 = E$ (in this case we say that $T$ has \emph{real multiplication}) or

 $\bullet$ $E_0 \subset E$ is a purely imaginary, quadratic extension and $'$ is the restriction of complex
  conjugation 
  to $E$ (we say that $T$ has \emph{complex multiplication}).
 \end{Thm}
\end{subsection}

\begin{subsection}{Splitting of Hodge structures of K3 type with real multiplication} 
(For this and
the next section see \cite{vG2}, 2.4 and 2.5.) \label{SplittingExtensions}
Let $(T,h,q)$ be an irreducible Hodge structure of K3 type and assume that 
$E = \End_{\Hdg}(T)$ is
a totally real number field. Note that 
by Theorem \ref{ThmZarhin1}, all endomorphisms in $E$ are $q$-self-adjoint.

By the theorem of the primitive element, there exists
$\alpha \in E$ such that $E = \QQ(\alpha)$. Let $d = [E: \QQ]$. 
Let $P$ be the minimal polynomial
of $\alpha$ over $\QQ$, denote by $\widetilde{E}$ the splitting field of $P$ in $\RR$.
Let $G = \mathrm{Gal}(\widetilde{E}/\QQ)$ and $H = \mathrm{Gal}(\widetilde{E}/E)$.
Choose $\sigma_1 = \id, \sigma_2, \ldots, \sigma_d \in G$ such that
\begin{equation*}
 G = \sigma_1 H \sqcup \ldots \sqcup \sigma_d H.
\end{equation*}
Note that each coset $\sigma_i H$ induces a well-defined embedding $E \hookrightarrow \widetilde{E}$.
In $\widetilde{E}[X]$ we get
\begin{equation*} 
 P(X) = \prod_{i=1}^d (X - \sigma_i(\alpha))
\end{equation*}
and consequently
\begin{equation*}
 \begin{aligned} E \otimes_{\QQ} \widetilde{E} = & \QQ[X]/(P) \otimes_{\QQ} \widetilde{E} \\ 
   \simeq & \bigoplus_{i=1}^d \widetilde{E}[X] / (X - \sigma_i (\alpha)) \\
   \simeq & \bigoplus_{i=1}^d E_{\sigma_i}.
 \end{aligned}
\end{equation*}
The symbol $E_{\sigma_i}$ stands for the field $\widetilde{E}$, the index $\sigma_i$ keeps track of the fact
that the $\widetilde{E}$-linear extension of $E \subset \End_{\QQ}(E)$ acts on $E_{\sigma_i}$ via $e(x) = \sigma_i(e) \cdot x$.
See Section \ref{SectionCorestriction} for another interpretation of $E_{\sigma_i}$.

In the same way, since $T$ is a finite-dimensional $E$-vector space we get a decomposition
\begin{equation*}
 T_{\widetilde{E}} = T \otimes_{\QQ} \widetilde{E}  = \bigoplus_{i=1}^d T_{\sigma_i}.
\end{equation*}
This is the decomposition of $T_{\widetilde{E}}$ into eigenspaces of the $\widetilde{E}$-linear
extension of the $E$-action on $T$, $T_{\sigma_i}$ being the eigenspace of $e_{\widetilde{E}}$ 
to the eigenvalue $\sigma_i (e)$ for $e \in E$. Since each $e \in E$ is $q$-self-adjoint (that is $e'= e$), 
the decomposition is orthogonal. Let $q_{\widetilde{E}}$ be the $\widetilde{E}$-bilinear extension of $q$ to $T_{\widetilde{E}} \times T_{\widetilde{E}}$.
Using the notation
\begin{equation*}
  T_i := T_{\sigma_i} \; \text{and} \; q_i = (q_{\widetilde{E}})_{| T_i \times T_i},
\end{equation*}
we have an orthogonal decomposition
\begin{equation} \label{ZerlegungT}
 (T_{\widetilde{E}}, q_{\widetilde{E}}) = \bigoplus_{i=1}^d (T_i, q_i).
\end{equation}
\end{subsection}

\begin{subsection}{Galois action on \texorpdfstring{$T_{\widetilde{E}}$}{T}} \label{GaloisT}
Letting $G$ act in the natural way on $\widetilde{E}$, we get a (only $\QQ$-linear) 
Galois action on $T_{\widetilde{E}} = T \otimes_{\QQ} \widetilde{E}$.
Under this action, for $\tau \in G$ we have 
\begin{equation} \label{GalPerm1}
 \tau T_{\sigma_i} = T_{ \tau \sigma_i}.
\end{equation}
This is because the Galois action commutes with the $\widetilde{E}$-linear extension of any endomorphism 
$e \in E \subset \End_{\QQ}(T)$ the latter being defined over $\QQ$ and because
for $t_i \in T_{\sigma_i}$ and $e \in E$
\begin{equation*}
  e_{\widetilde{E}} (\tau (t_i)) = 
  \tau (e_{\widetilde{E}} (t_i)) = \tau (\sigma_i(e) t_i) = \tau(\sigma_i (e)) \tau(t_i)
   = (\tau \sigma_i (e)) \tau (t_i),
\end{equation*}
which means that $\tau$ permutes the eigenspaces of $e_{\widetilde{E}}$ precisely in the way we claimed.
Define a homomorphism
\begin{equation} \label{Perm1-d} \gamma: G \to \; \mathfrak{S}_d, \; \; \;
 \tau \mapsto \; \{ i \mapsto \tau (i) \; \text{where} \; (\tau \sigma_i) H = \sigma_{\tau(i)} H \}.
\end{equation}
(This describes the action of $G$ on $G/H$.)
With this notation, (\ref{GalPerm1}) reads 
\begin{equation} \label{GalPerm}
  \tau T_i = T_{\tau(i)}.
\end{equation}

Interpret $T$ as a subspace of $T_{\widetilde{E}}$ via the natural inclusion 
$T \hookrightarrow T_{\widetilde{E}}, \; t \mapsto t \otimes 1$. Denote by $\pi_i$ the projection to $T_i$.
For $t \in T$ and $\tau \in G$ we have $t = \tau (t)$. Write $t_i := \pi_i (t \otimes 1)$, then
$t = \sum_i t_i$. Using (\ref{GalPerm})
we see that
\begin{equation} \label{t_i=sigma_i}
 t_{\tau i} = \tau (t_i).
\end{equation}
It follows that 
\begin{equation} \label{iota}
 \iota_i: T \to \; T_i, \; \; \; t \mapsto  \; \pi_i (t \otimes 1)
\end{equation}
is an injective map of $E$-vector spaces ($E$ acting on $T_i$ via $\sigma_i : E \hookrightarrow \widetilde{E}$).
Equation (\ref{t_i=sigma_i}) can be rephrased as
\begin{equation} \label{iota_i=sigma_i}
 \iota_{\tau i} = \tau \circ \iota_i.
\end{equation}

Since $q$ is defined over $\QQ$, we have for $t \in T_{\widetilde{E}}$ and $\tau \in G$
\begin{equation*}
 q_{\widetilde{E}} (\tau t) = \tau q_{\widetilde{E}}(t).
\end{equation*}
This implies that for $t \in T$
\begin{equation} \label{QuadratischeFormGalois}
 q_i (\iota_i(t)) = \sigma_i q_1(\iota_1(t)).
\end{equation}
\end{subsection}
\begin{subsection}{The special Mumford--Tate group of a Hodge structure of K3 type with real multiplication}
\label{SectionSMTZarhin}
  Zarhin \cite{Zarhin} also computed the special Mumford--Tate group of an irreducible Hodge
  structure of K3 type. Recall that for a Hodge structure $(W, h: \mathrm{U}(1) \to \mathrm{GL}(W)_{\RR})$ the
  special Mumford--Tate group $\mathrm{SMT}(W)$ is the smallest linear algebraic subgroup of $\mathrm{GL}(W)$ 
  defined over $\QQ$ with $h(\mathrm{U}(1)) \subset \mathrm{SMT}(W)_{\RR}$ (cf.\ \cite{Go}). 

  Assume that $(T,h,q)$ is an irreducible Hodge structure of K3 type with real multiplication by $E = \End_{\Hdg}(T)$.
  We continue to use the notations of Section \ref{SplittingExtensions}.
 Denote by $Q$ the restriction of $q_1$ to $T \subset T_1$ (use the inclusion $\iota_1$ of (\ref{iota})). This is an $E$-valued (since $H$-invariant), 
non-degenerate, symmetric bilinear form
on the $E$-vector space $T$. Denote by $\SO(Q)$ the $E$-linear algebraic group of 
$Q$-orthogonal, $E$-linear transformations of $T$ with determinant $1$. 

Recall that for
an $E$-variety $Y$ the \emph{Weil restriction} $\Res_{E /\QQ}(Y)$ is the $\QQ$-variety whose $K$-rational
points are the $E \otimes_{\QQ} K$-rational points of $Y$ for any extension field $\QQ \subset K$ (cf.\ \cite{BLR}).

\begin{Thm}[Zarhin, see \cite{Zarhin}, see also \cite{vG2}, 2.8] \label{Thm2Zarhin}
The special Mumford--Tate group of the Hodge structure $(T,h,q)$ with real multiplication by $E$ is \begin{equation*}
  \SMT(T) = \Res_{E/\QQ}(\SO(Q)). 
\end{equation*}
Its representation on $T$ is the natural one, where we regard $T$ as
a $\QQ$-vector space and use that any $E$-linear endomorphism of $T$ is in particular $\QQ$-linear.
After base change to $\widetilde{E}$
\begin{equation*}
  \SMT(T)_{\widetilde{E}} = \prod_i \SO((T_i), (q_i)),
\end{equation*}
its representation on $T_{\widetilde{E}} = \bigoplus_i (T_i)$ is the product 
of the standard representations. 
\end{Thm}
\end{subsection}
\end{section}

\begin{section}{Kuga--Satake varieties and real multiplication} \label{SectionKS}
\begin{subsection}{Kuga--Satake varieties}
Let $(T,h,q)$ be a Hodge structure of K3 type. Kuga and Satake \cite{KS} found
a way to associate to
this a polarizable $\QQ$-Hodge structure of weight one $(V, h_s : \mathrm{U}(1) \to \mathrm{GL}(V)_{\RR})$, 
in other words an isogeny class of Abelian varieties, together with an inclusion of Hodge structures 
\begin{equation} \label{EmbeddingTransKS}
  T \subset V \otimes V.
\end{equation} 
Set $V:= C^0(q)$ where $C^0(q)$ is the even Clifford algebra of $q$.
Define a weight one Hodge structure on $V$ in the following way: 
Choose $f_1, f_2 \in (T^{2,0} \oplus T^{0,2})_{\RR}$ such that $\CC (f_1 + i f_2) = T^{2,0}$ and $q(f_i, f_j) =  \delta_{i,j}$ (recall that
$q_{|(T^{2,0} \oplus T^{0,2})_{\RR}}$ is positive definite). Define $J: V \to V, \; v \mapsto f_1 f_2 v$, then we see that $J^2 = - \id$. 
Now we can define a homomorphism of algebraic groups
\begin{equation*}
  h_s: \UU(1) \to \mathrm{GL}(V)_{\RR}, \; \; \exp(x i) \mapsto \exp(x J),
\end{equation*}
and this induces the Kuga--Satake Hodge structure. One can check that $h_s$ is independent of the choice of $f_1, f_2$
(see \cite[Lemma 5.5]{vG1}).

It can be shown that the Kuga--Satake Hodge structure admits a polarization (cf.\ \cite[Prop.\ 5.9]{vG1}) and
that there is an embedding of Hodge structures as in (\ref{EmbeddingTransKS}) (see \cite[Prop.\ 6.3]{vG1}).
\end{subsection}

\begin{subsection}{Corestriction of algebras}  \label{SectionCorestriction}
Let $E/K$ be a finite, separable extension of fields 
of degree $d$ and
let $A$ be an $E$-algebra. We use the notations of Section \ref{SplittingExtensions}, 
so $\widetilde{E}$ is a 
normal closure of $E$ over $K$, $\sigma_1, \ldots ,\sigma_d$ is a
set of representatives of $G/H$ where $G = \mathrm{Gal}(\widetilde{E}/K)$ and $H = \mathrm{Gal}(\widetilde{E}/E)$.

For $\sigma \in G$ 
define the twisted $\widetilde{E}$-algebra as the ring
\begin{equation*}
  A_{\sigma} := A \otimes_E \widetilde{E}
\end{equation*}
which carries an $\widetilde{E}$-algebra structure given by
\begin{equation*}
  \lambda \cdot (a \otimes e) = a \otimes \sigma^{-1} (\lambda) e.
\end{equation*}
Note that $A_{\sigma} \simeq A \otimes_E E_{\sigma}$.

Let $V$ be an $E$-vector space and $W$ an $\widetilde{E}$-vector space,
let $\sigma \in G$.
A homomorphism of $K$-vector spaces $\varphi: V \to W$ is called \emph{$\sigma$-linear}
if $\varphi (\lambda v) = \sigma (\lambda) \varphi (v)$ for all $v \in V$ and $\lambda \in E$.
If both, $V$ and $W$ are $\widetilde{E}$-vector spaces, there is a similar notion of an $\sigma$-linear
homomorphism.

\begin{Lemma} \label{universal_twisted}
  The map
  \begin{equation*}
    \kappa_{\sigma}: A \to  A_{\sigma}, \; \; a \mapsto a \otimes 1
  \end{equation*}
  is a $\sigma$-linear ring homomorphism and the pair $(A_{\sigma}, \kappa_{\sigma})$ has the
  following universal property: For all $\widetilde{E}$-algebras
  $B$ and for all $\sigma$-linear ring
  homomorphisms $\varphi: A \to B$ there exists a unique $\widetilde{E}$-algebra
  homomorphism $\widetilde{\varphi}: A_{\sigma} \to B$ making the diagram
  \begin{equation*}
    \xymatrix{ A \ar[r]^-{\kappa_{\sigma}} \ar[dr]_{\varphi} &  A_{\sigma} \ar[d]^{\widetilde{\varphi}} \\
               & B }
  \end{equation*}
  commutative.
\end{Lemma}

\begin{Proof} We only check the universal
property.
To give a $K$-linear map $\alpha : A \otimes_E \widetilde{E} \to B$
is the same as to give a $K$-bilinear map $\beta: A \times \widetilde{E} \to B$
satisfying 
\begin{equation} \label{beta_lambda}
\beta(\lambda a , e) = \beta(a ,\lambda e)
\end{equation}
for all $a \in A, e \in \widetilde{E}$ and $\lambda \in E$.
These maps are related by the condition
\begin{equation*}
  \alpha (a \otimes e) = \beta(a,e).
\end{equation*}

Now given $\varphi$ as in the lemma, we define
  \begin{equation*}
   \psi: A \times \widetilde{E} \to B, \; \; (a,e) \mapsto \sigma(e) \varphi(a).
  \end{equation*}
This is a $K$-bilinear map which satisfies (\ref{beta_lambda}) and therefore,
it induces a $K$-linear map
\begin{equation*}
  \widetilde{\varphi}: A \otimes_E \widetilde{E} \to B, \; \; a \otimes e \mapsto \sigma(e) \varphi(a).
\end{equation*}
It is clear that $\widetilde{\varphi}$ is a ring homomorphism and that it respects
the $\widetilde{E}$-algebra structures if we interpret $\widetilde{\varphi}$ as a map
$\widetilde{\varphi}:  A_{\sigma} \to B$. The uniqueness of this map is immediate.
\end{Proof}

\vspace{2ex} \noindent
{\it Remark.} (i) The lemma shows that up to unique $\widetilde{E}$-algebra isomorphism,
the twisted algebra $ A_{\sigma_i}$ depends only
on the coset $\sigma_i H$. Indeed, for $\sigma \in \sigma_i H $ the inclusion 
$A \hookrightarrow  A_{\sigma_i} $ is $\sigma$-linear
because $\sigma$ and $\sigma_i$ induce the same
embedding of $E$ into $\widetilde{E}$. By the lemma, there exists an $\widetilde{E}$-algebra
isomorphism $\alpha_{\sigma, \sigma_i} : A_{\sigma} \stackrel{\sim}{\to}  A_{\sigma_i},
\; a \otimes e \mapsto \sigma(e) \cdot (a \otimes 1) = a \otimes \sigma_i^{-1} \sigma (e)$.

\vspace{1ex}
(ii) In Section \ref{SplittingExtensions} we were in the situation $E = \QQ(\alpha)$. There
we discussed the splitting $E \otimes_{\QQ} \widetilde{E} \simeq \bigoplus_i 
\widetilde{E}[X] / (X - \sigma_i(\alpha)) \simeq \bigoplus_i E_{\sigma_i}$ and we used the symbol 
$E_{\sigma_i}$ for the field $\widetilde{E}$ with $E$-action via $e(x) = \sigma_i (e) \cdot x$. This is
precisely our twisted $\widetilde{E}$-algebra $E_{\sigma_i}$ on which $E$ acts via the inclusion $\kappa_{\sigma_i}$.

\vspace{2ex}
For $\tau \in G$ there is a unique $\tau$-linear ring
isomorphism $\tau:  A_{\sigma_i} \to  A_{\sigma_{\tau i}}$ which extends the identity on 
$A \subset  A_{\sigma_i}$ (in the sense that
$\tau \circ \kappa_{\sigma_i} = \kappa_{\sigma_{\tau i}}$). This map is given as the composition of the following two maps:
First apply the identity map $ A_{\sigma_i} \to  A_{\tau \sigma_i}, \; a \otimes e \mapsto a \otimes e$ which is 
a $\tau$-linear ring isomorphism. Then apply the isomorphism $\alpha_{\tau \sigma_i, \sigma_{\tau i}}$ (by
definition of the $G$-action on $\{1, \ldots, d\}$ we have $\tau \sigma_i \in \sigma_{\tau i} H$).
On simple tensors the map $\tau$ 
takes the form
\begin{equation} \label{taulinear}
 a \otimes e \mapsto a \otimes \sigma_{\tau i}^{-1} \tau \sigma_i (e).
\end{equation}

These maps induce a natural action of $G$ on
\begin{equation*}
 Z_G(A) :=  A_{\sigma_1} \otimes_{\widetilde{E}} \ldots \otimes_{\widetilde{E}} A_{\sigma_d}
\end{equation*}
where
\begin{multline} \label{GWirkungCores}
  \tau ((a_1 \otimes e_1) \otimes \ldots \otimes (a_d \otimes e_d))  \\
   = \big( a_{\tau^{-1}1} \otimes \sigma_1^{-1} \tau \sigma_{\tau^{-1}1} (e_{\tau^{-1}1}) \big) \otimes \ldots
   \otimes \big( a_{\tau^{-1}d} \otimes \sigma_d^{-1} \tau \sigma_{\tau^{-1}d} (e_{\tau^{-1} d}) \big).
\end{multline}

\begin{Def}[\cite{D}, \S 8, Def.\ 2 or \cite{T}, 2.2] The \emph{corestriction of $A$ to $K$}
 is the $K$-algebra of $G$-invariants in $Z_G(A)$
\begin{equation*}
 \cores_{E/K} (A): = Z_G(A)^G.
\end{equation*}
\end{Def}

\noindent {\it Remark.} 
(i) By \cite[\S 8, Cor.\ 1]{D} there is a natural isomorphism
\begin{equation*}
  \cores_{E/K} (A) \otimes_{K} \widetilde{E} \simeq Z_G(A)
\end{equation*}
In particular, with $d = [E: K]$ one gets $\dim_K \cores_{E/K} (A) = (\dim_E(A))^d$.

\vspace{1ex}
(ii) Let $X = \mathrm{Spec}(A)$ for a commutative $L$-algebra $A$. Then for any $K$-algebra $B$ we get a chain of
 isomorphisms, functorial in $B$
 \begin{equation*} \begin{aligned}
   \Hom_{K-\mathrm{Alg}} (\cores_{E/K}(A), B) & \simeq  \; 
   \left( \Hom_{\widetilde{E}-\mathrm{Alg}} (Z_G(A), B \otimes_K \widetilde{E}) \right)^G  \\
   & \simeq \;
   \Hom_{E-\mathrm{Alg}}(A, B \otimes_K E). \end{aligned}
 \end{equation*}
 Here, the last isomorphism is given by composing $f \in \left( \Hom_{\widetilde{E}-\mathrm{Alg}} (Z_G(A), B \otimes_K \widetilde{E}) \right)^G$
 with the inclusion $j: A \hookrightarrow Z_G(A), \; a \mapsto \kappa_{\sigma_1}(a) \otimes 1 \otimes \ldots 
 \otimes 1$. 
 (The image of this composition is contained in the $H$-invariant part of $B \otimes_K \widetilde{E}$ which is $B \otimes_K E$.)
 This map is an isomorphism, since $Z_G(A)$ is generated as an $\widetilde{E}$-algebra by elements of the form
 $\sigma \circ j (a)$ with $a \in A$ and $\sigma \in G$.

 It follows that
 \begin{equation*}
   \mathrm{Res}_{E/K} (\mathrm{Spec}(A)) \simeq \mathrm{Spec}(\cores_{E/K} (A)),
 \end{equation*}
 i.e.\ the Weil restriction of affine $E$-schemes is the same as the corestriction of commutative $E$-algebras.

\end{subsection}

\begin{subsection}{The decomposition theorem}
We will assume from now to the end of the section that $(T,h,q)$ is an irreducible Hodge structure
of K3 type with $E=\End_{\Hdg}(T)$ a totally real number field.

Recall that in this case $T$ is an $E$-vector space which carries a natural $E$-valued quadratic form
$Q$ (see \ref{SectionSMTZarhin}). Let $C^0(Q)$ be the even Clifford algebra of $Q$ over $E$. It was van Geemen (see \cite[Prop.\ 6.3]{vG2}) who discovered that the algebra $\cores_{E/\QQ} (C^0(Q))$
appears as a sub-Hodge structure in the Kuga--Satake Hodge structure of $(T,h,q)$. We are going to show that
this contains all information on the Kuga--Satake Hodge structure.

\begin{Thm} \label{Satz} 
 Denote by $(V,h_s)$ the Kuga--Satake Hodge structure of $(T,h,q)$.
 
 \vspace{1ex}
 \emph{(i)} The special Mumford--Tate group of $(V,h_s)$ is the image of $\Res_{E/\QQ} (\Spin(Q))$
 in $\Spin(q)$ under a morphism $m$ of rational algebraic groups which after 
 base change to $\widetilde{E}$ becomes
 \begin{equation*} m_{\widetilde{E}}:
   \Spin(q_1) \times \ldots \times \Spin(q_d) \to  \Spin(q)_{\widetilde{E}}, \; \; \;
   (v_1, \ldots, v_d) \mapsto v_1 \cdot \ldots \cdot v_d. 
 \end{equation*}
 \vspace{1ex}
 \emph{(ii)} Let $W := \cores_{E/\QQ}(C^0(Q))$. Then $W$ can be canonically embedded in $V$ and
 the image is $\SMT(V)$-stable and therefore, it is a sub-Hodge structure. 
 Furthermore, there is a (non-canonical) isomorphism of Hodge
 structures
 \begin{equation*}
   V \simeq W^{2^{d-1}}.
 \end{equation*}
 
 \vspace{1ex}
 \emph{(iii)} We have
 \begin{equation*}
  \End_{\Hdg}(W) = \cores_{E/\QQ}(C^0(Q))
 \end{equation*}
 and consequently
 \begin{equation*}
  \End_{\Hdg}(V) = \mathrm{Mat}_{2^{d-1}} \big( \cores_{E/\QQ} (C^0(Q)) \big).
 \end{equation*} 
\end{Thm}
The proof will be given in Section \ref{proof}. The theorem tells us that the Kuga--Satake
variety $A$ of $(T,h,q)$ is isogenous to a self-product $B^{2^{d-1}}$ of an Abelian
variety $B$ with $\End_{\QQ}(B) = \cores_{E/\QQ} (C^0(Q))$ and therefore, it proves Theorem
\ref{ThmZerlegungKS}. 

Note that $B$ is not 
simple in general. We will see examples below where $B$ decomposes further.
\end{subsection}

\begin{subsection}{Galois action on \texorpdfstring{$C(q)_{\widetilde{E}}$}{C(q)}}
By Section \ref{SplittingExtensions} there is a decomposition
\begin{equation*}
 (T,q)_{\widetilde{E}} = \bigoplus_{i=1}^d (T_i, q_i).
\end{equation*}
This in turn yields an isomorphism
\begin{equation*}
  C(q)_{\widetilde{E}} \simeq C (q_{\widetilde{E}})    \simeq  
  C(q_1) \widehat{\otimes}_{\widetilde{E}} \ldots \widehat{\otimes}_{\widetilde{E}} C(q_d).
\end{equation*}
Here, the symbol $\widehat{\otimes}$ denotes the graded tensor product of algebras, which on the
level of vector spaces is just the usual tensor product, but which twists the algebra structure by a 
suitable sign (see \cite[??]{vG1}).

Decompose $C(q_i) = C^0(q_i) \oplus C^1(q_i)$ in the even and the odd part.
If we forget the algebra structure and only look at $\widetilde{E}$-vector spaces, we get
\begin{equation*}
  C(q)_{\widetilde{E}} = 
  \bigoplus_{\mathbf{a} \in \{0,1 \}^d} C^{a_1} (q_1) \otimes_{\widetilde{E}} \ldots \otimes_{\widetilde{E}} C^{a_d} (q_d).
\end{equation*}
For $\mathbf{a} = (a_1,\ldots, a_d) \in \{0,1 \}^d$ define
\begin{equation*}
 C^{\mathbf{a}}(q) = C^{a_1}(q_1) \otimes \ldots \otimes C^{a_d}(q_d).
\end{equation*}
We introduced an action of $G = \mathrm{Gal}(\widetilde{E}/\QQ)$ on $\{1, \ldots , d\}$ 
(see (\ref{Perm1-d})). This induces an action
\begin{equation} \label{WirkungG0,1^d}
 G \times \{0,1 \}^d \to \{ 0,1 \}^d, \; \; (\tau, (a_1, \ldots, a_d)) \mapsto (a_{\tau^{-1} 1}, \ldots, a_{\tau^{-1} d}).
\end{equation} 
The next lemma describes the Galois action on $C(q)_{\widetilde{E}}$.
\begin{Lemma} \label{GaloisClifford}
  \emph{(i)} Via the map
  \begin{equation*} 
     C(q_i) \subset C(q)_{\widetilde{E}}, \; v_i \mapsto 1 \otimes \ldots \otimes v_i \otimes
  \ldots \otimes 1
  \end{equation*}
  we interpret $C(q_i)$ as a subalgebra of $C(q)_{\widetilde{E}}$. Then the restriction
  of $\tau \in G$ to $C(q_i)$ induces
  an isomorphism of $\ZZ / 2\ZZ$-graded $\QQ$-algebras
  \begin{equation*}
   \tau: (C(q_i)) \stackrel{\sim}{\to} C(q_{\tau{i}}).
  \end{equation*}
 \vspace{1ex}
 \emph{(ii)} For $\tau \in G$ and $\mathbf{a} \in \{ 0,1 \}^d$ we get
  \begin{equation*}
    \tau (C^{\mathbf{a}}(q)) = C^{\tau \mathbf{a}} (q).
  \end{equation*}
\end{Lemma}

\begin{Proof} 
 Tensor the natural inclusion $T \hookrightarrow C(q)$ with $\widetilde{E}$ to get a
 $G$-equivariant inclusion
 \begin{equation*}
  T \otimes_{\QQ} \widetilde{E} = \bigoplus_{i=1}^d T_i \to C(q)_{\widetilde{E}}.
 \end{equation*}
 Using (\ref{GalPerm}), we find for $t_i \in C(q_i)$ that 
 $\tau (t_i) \in C(q_{\tau(i)})$. Now, $C(q_i)$ is spanned as a $\QQ$-algebra
 by products of the form 
\begin{equation*}
  t_1 \cdot \ldots \cdot t_k
  = \pm (1 \otimes \ldots \otimes t_1 \otimes \ldots \otimes 1) \cdot \ldots \cdot (1 \otimes \ldots \otimes t_k \otimes
\ldots \otimes 1)
\end{equation*}
for $t_1, \ldots, t_k \in T_i$. Since $G$ acts by $\QQ$-algebra homomorphisms on $C(q)_{\widetilde{E}}$, this implies (i). 

Item (ii) is an immediate consequence of (i): The space $C^{\mathbf{a}}(q)$ is spanned as $\QQ$-vector space
by products of the form 
$v_1 \cdot \ldots \cdot v_d = \pm v_1 \otimes \ldots \otimes v_d$ with $v_i \in C^{a_i}(q_i)$. Then use again, that $G$ acts by $\QQ$-algebra homomorphisms.
\end{Proof}

\begin{Lemma} \label{twistedC(Q)}
 For $i \in \{1, \ldots, d \}$ the twisted algebra $C^0(Q)_{\sigma_i}$ is canonically isomorphic 
 as an $\widetilde{E}$-algebra to
 $C^0(q_i)$. Thus 
 \begin{equation*}
   Z_G(C^0(Q)) \simeq C^0(q_1) \otimes_{\widetilde{E}} \ldots \otimes_{\widetilde{E}} C^0(q_d).
 \end{equation*}
 On both sides there are natural $G$-actions: On the left hand side $G$ acts via the action introduced in 
 (\ref{GWirkungCores}), whereas on the right hand side it acts via the restriction of its action on
 $C(q)_{\widetilde{E}}$ (use Lemma \ref{GaloisClifford}). Then the above isomorphism is $G$-equivariant.
\end{Lemma}

\begin{Proof} Fix $i \in \{1, \ldots, d\}$.
 The composition of the canonical inclusion
 $C^0(Q) \subset C^0(q_1) \simeq C^0(Q)_{\widetilde{E}}$ with the
 restriction to $C^0(Q)$ of the map $\sigma_i : C(q_1) \to C(q_i)$
 from Lemma \ref{GaloisClifford} induces a $\sigma_i$-linear ring homomorphism
 \begin{equation*}
   \varphi_i : C^0(Q) \hookrightarrow C^0(q_i).
 \end{equation*}
 By Lemma \ref{universal_twisted} we get an $\widetilde{E}$-algebra homomorphism
 \begin{equation*}
   \widetilde{\varphi}_i : C^0(Q)_{\sigma_i} \to C^0(q_i).
 \end{equation*}
 Recall that there are inclusions $\iota_i: T \hookrightarrow T_i$ (see (\ref{iota}))
 which satisfy $\tau \circ \iota_i = \iota_{\tau i}$ (see (\ref{iota_i=sigma_i})). 
 Let $t_1, \ldots, t_m \in T$ such that $\iota_1(t_1), \ldots, \iota_1(t_m)$ form
 a $q_1$-orthogonal basis of $T_1$. Then the vectors
 $\iota_i (t_1), \ldots, \iota_i(t_m)$ form a $q_i$-orthogonal basis of $T_i$ 
 (use (\ref{QuadratischeFormGalois})). By definition of
 $\widetilde{\varphi_i}$
 \begin{equation} \label{WirkungTildePhi}
   \widetilde{\varphi}_i \left( \iota_1 (t_1)^{i_1} \cdot \ldots \cdot \iota_1 (t_m)^{i_d} \right) =
    \iota_i(t_1)^{i_1} \cdot \ldots \cdot \iota_i (t_m)^{i_m}.
 \end{equation}
 This implies that $\widetilde{\varphi}_i$ maps an $\widetilde{E}$-basis of $C^0(Q)_{\sigma_i}$
 onto an $\widetilde{E}$-basis of $C^0(q_i)$, whence $\widetilde{\varphi}_i$ is an isomorphism
 of $\widetilde{E}$-algebras.
 
 As for the $G$-equivariance, we have to check that for all $\tau \in G$ the diagram
 \begin{equation*}
 \begin{CD} 
    C^0(Q)_{\sigma_i} @>\widetilde{\varphi}_i>> C^0(q_i) \\
    @V{\tau}VV @VV{\tau}V \\
    C^0(Q)_{\sigma_{\tau i}} @>\widetilde{\varphi}_{\tau i}>> C^0(q_{\tau i})
 \end{CD}
 \end{equation*}
 is commutative. It is enough to check this on an $\widetilde{E}$-basis of $C^0(Q)_{\sigma_i}$ because the
 vertical maps are both $\tau$-linear whereas the horizontal ones are $\widetilde{E}$-linear. 
 Since $\tau: C^0(Q)_{\sigma_i} \to C^0(Q)_{\sigma_{\tau i}}$ was defined as the extension of
 the identity map on $C^0(Q) \subset C^0(Q)_{\sigma_i}$, we have
 \begin{equation*} \begin{aligned}
   \widetilde{\varphi}_{\tau i} \circ \tau \left( \iota_1 (t_1)^{i_1} \cdot \ldots \cdot \iota_1(t_m)^{i_m} \right) = & 
    \widetilde{\varphi}_{\tau i} \left( \iota_1 (t_1)^{i_1} \cdot \ldots \cdot \iota_1(t_m)^{i_m} \right) \\
    = & \; \iota_{\tau i}(t_1)^{i_1} \cdot \ldots \cdot \iota_{\tau i} (t_m)^{i_m} \\
    = & \; (\tau \circ \iota_i)(t_1)^{i_1} \cdot \ldots \cdot (\tau \circ \iota_i) (t_m)^{i_m} \\
    = & \; \tau \left( \iota_i (t_1)^{i_1} \cdot \ldots \cdot \iota_i (t_m)^{i_m} \right) \\
    = & \; \tau \circ \widetilde{\varphi}_i \left( \iota_1 (t_1)^{i_1} \cdot \ldots \cdot \iota_m (t_m)^{i_m} \right).
 \end{aligned} \end{equation*}
 This completes the proof of the lemma. 
\end{Proof}
\end{subsection}

\begin{subsection}{Proof of the decomposition theorem} \label{proof}
Let $K$ be a field and $(U,r)$ be a quadratic $K$-vector space.
Recall that the spin group of $r$ comes with two natural representations:

First there is the covering representation $\rho: \Spin(r) \to \SO(r)$ which over 
an extension field $K \subset L$ maps $y \in \Spin(r)(L) = \{ x \in (C^0(r) \otimes_{K} L)^* \; | \; x \iota(x) = 1 \; \text{and} \;
x U x^{-1} \subset U\}$ to the endomorphism $U \to U, \; u \mapsto x u x^{-1}$. Here, $\iota : C(r) \to C(r)$
is the natural involution of the Clifford algebra.

Secondly, the spin
representation realizes $\Spin(r)$ as a subgroup of $\mathrm{GL}(C^0(r))$ by sending
$y \in \Spin(r)(L)$ to the endomorphism of $C^0(r)$ given by $x \mapsto y \cdot x$.

\vspace{2ex}

\noindent \it Proof of (i). \rm By \cite[Prop.\ 6.3]{vG1}, there is a commutative diagram
\begin{equation} \label{DiagrammHodge}
  \begin{CD} \UU(1) @>h_s>> \Spin(q)_{\RR} @>>> \mathrm{GL}(C^0(q))_{\RR} \\
   @| @VV{\rho}V \\
   \UU(1) @>>{h}> \SO(q)_{\RR} @>>> \mathrm{GL}(T)_{\RR}.
  \end{CD}
\end{equation}
(Van Geemen works with the Mumford--Tate group, therefore he gets a factor $t^2$ in 6.3.2. 
This factor is 1 if one restricts the attention to the special Mumford--Tate group, moreover 
it is then clear that $h_s(\CC^*) \subset \mathrm{CSpin(q)} = \{ v \in C^0(q)^* \; | \; v T v^{-1} \subset T \}$ 
implies $h_s(\UU(1)) \subset \Spin(q)$.)

\vspace{2ex}

\noindent \bf Claim: \rm There is a Cartesian diagram
\begin{equation*}
 \begin{CD} \SMT(V) @>>> \Spin(q) \\
            @V{\rho_{|\SMT(V)}}VV @VV{\rho}V \\
            \SMT(T) @>>> \SO(q).
 \end{CD}
\end{equation*}
where the horizontal maps are appropriate factorizations of the inclusions $\SMT \subset \mathrm{GL}$ whose
existence is guaranteed by (\ref{DiagrammHodge}).
\vspace{2ex}

\noindent \it Proof of the claim. \rm It is clear by looking at (\ref{DiagrammHodge}) and at
the definition of the special Mumford--Tate group that 
\begin{equation*}
 \SMT(V) \subset \SMT(T) \times_{\SO(q)} \Spin(q).
\end{equation*}
In the same way we see that
\begin{equation*}
 \SMT(T) \subset \rho(\SMT(V))
\end{equation*}
and hence we have a chain of inclusions
\begin{equation*}
 \SMT(V) \subset \SMT(T) \times_{\SO(q)} \Spin(q) \subset \rho(\SMT(V)) \times_{\SO(q)} \Spin(q).
\end{equation*}
But over any field, the kernel of $\rho$ consists of $\{ \pm 1 \} \subset \SMT(V)$ (because $h_s(-1) = -1$)
and thus
\begin{equation*}
 \SMT(V) = \rho(\SMT(V)) \times_{\SO(q)} \Spin(q).
\end{equation*}
This proves the claim. \hfill $\mathrm{(Claim)} \Box$  \par

\vspace{2ex}
To continue the proof of (i) we have to define the morphism of rational algebraic groups
\begin{equation*}
 m: \Res_{E/\QQ}(\Spin(Q)) \to \Spin(q).
\end{equation*}
For that sake, note first that there is a natural isomorphism of $\widetilde{E}$-algebras
\begin{equation} \label{C(Q)ext} \begin{aligned}
  C^0(Q) \otimes_{\QQ} \widetilde{E} \simeq & \;C^0(Q) \otimes_E (E \otimes_{\QQ} \widetilde{E}) \\
  \simeq & \; \bigoplus_i C^0(Q) \otimes_E E_{\sigma_i} \\
  \simeq  & \; \bigoplus_i C^0(Q)_{\sigma_i} \\ \simeq & \;  C^0(q_1) \oplus \ldots \oplus C^0(q_d) \end{aligned}
\end{equation}
where we use the notations of Section \ref{SectionCorestriction} and for the last identification Lemma \ref{twistedC(Q)}.
Consider the natural $G$-action on $C^0(q_1) \oplus \ldots \oplus C^0(q_d)$ given by
\begin{equation*} 
 (\tau, (v_1, \ldots, v_d)) \mapsto (\tau v_{\tau^{-1} 1} , \ldots , \tau v_{\tau^{-1}d}).
\end{equation*} 
On $C^0(Q)\otimes_{\QQ} \widetilde{E}$, the Galois group $G$ acts by its natural action on $\widetilde{E}$. Then
the identification made in (\ref{C(Q)ext}) is $G$-equivariant and
we get an isomorphism of $\QQ$-vector spaces
\begin{equation*} C^0(Q) \simeq \left( C^0(q_1) \oplus \ldots \oplus C^0(q_d) \right)^G, \; \; \;
 v \mapsto (\sigma_1(v), \ldots, \sigma_d(v)). 
\end{equation*}
Now, look at the morphism of $\widetilde{E}$-affine spaces 
\begin{equation*} C^0(q_1) \oplus \ldots \oplus C^0(q_d) \to C^0(q)_{\widetilde{E}}, \; \; \; (v_1, \ldots, v_d) \mapsto v_1 \cdot \ldots  \cdot v_d.  
\end{equation*}
This morphism is $G$-equivariant on the $\widetilde{E}$-points and hence it comes from a
morphism of $\QQ$-varieties
\begin{equation*}
 \Res_{E/\QQ} C^0(Q) \to C^0(q).
\end{equation*}
The restriction of this latter to $\Res_{E/\QQ} (\Spin(Q))$ is the morphism $m$ we are
looking for. It is a morphism of algebraic groups which
after base change to $\widetilde{E}$ takes the form
\begin{equation*}
 m_{\widetilde{E}}: \Res_{E/\QQ}(\Spin(Q))_{\widetilde{E}} \simeq \Spin(q_1) \times \ldots \times \Spin(q_d)
  \to \Spin(q)_{\widetilde{E}}, \; \; \; (v_1, \ldots, v_d) \mapsto v_1 \cdot \ldots \cdot v_d.
\end{equation*}
It remains to show that the image of $m$ in $\Spin(q)$ is $\SMT(V)$. Using the
claim we have to show that the following diagram exists and that it is Cartesian
\begin{equation} \label{DiagrammSMT}
 \begin{CD}
  \im (m) @>>> \Spin(q) \\
  @V{\rho_{| \im(m)}}VV @VV{\rho}V \\
  \Res_{E / \QQ}(\SO(Q)) @>>> \SO(q).
 \end{CD}
\end{equation}
Here, the lower horizontal map is the one coming from Zarhin's Theorem \ref{Thm2Zarhin}.  

 It is enough to study (\ref{DiagrammSMT}) on $\overline{\QQ}$-points. It is easily seen that over $\widetilde{E} \subset \overline{\QQ}$
 the composition $\rho \circ m$ factorizes over 
 \begin{equation*}
   \rho_1 \times \ldots \times \rho_d : \Spin(q_1) \times \ldots \times \Spin(q_d) \to \SO(q_1) \times \ldots \times \SO(q_d) \simeq \Res_{E/\QQ}(\SO(Q))_{\widetilde{E}} \subset \SO(q)_{\widetilde{E}}.
 \end{equation*} 
 This shows that (\ref{DiagrammSMT}) exists. Moreover we see that $\rho_{|\im(m)}$ surjects onto $\SMT(T)(\overline{\QQ})$
 because $\rho_1 \times \ldots \times \rho_d$ does so. Since $\ker(\rho) = \{ \pm 1 \} \subset \im(m)$,
 the diagram (\ref{DiagrammSMT}) is Cartesian. This completes the proof of (i).
 \hfill    $\mathrm{(i)} \Box$  \par

\vspace{2ex}
\noindent \it Proof of (ii). \rm Choose $\mathbf{a}_0 = (0, \ldots, 0), \ldots, \mathbf{a}_r \in \{0,1 \}^d$ such
that 
\begin{equation*}
 \left\{ \mathbf{a} \in \{ 0,1\}^d \; | \; \sum_i a_i \equiv 0 \; (2) \right\} = 
  G \mathbf{a}_0 \sqcup \ldots \sqcup G \mathbf{a}_r,
\end{equation*}
where $G$ acts on $\{ 0,1 \}^d$ via the action introduced in (\ref{WirkungG0,1^d}).
Let $G_{\mathbf{a}_j} \subset G$ be the stabilizer of $\mathbf{a}_j$. Then 
\begin{equation} \label{C0=PlusDa} \begin{aligned}
  C^0(q)_{\widetilde{E}} = & \; \bigoplus_{j=0}^r \left( \bigoplus_{[\tau] \in G / G_{\mathbf{a}_j}} C^{\tau \mathbf{a}_j}(q)
   \right) \\
  = & \; \bigoplus_{j=0}^r D^{\mathbf{a}_j} \end{aligned}
\end{equation}
with $D^{\mathbf{a}_j} = \bigoplus_{[\tau] \in G / G_{\mathbf{a}_j}} C^{\tau \mathbf{a}_j}(q)$. 

By Lemma 
\ref{GaloisClifford} this is a decomposition of $G$-modules. Moreover, recall that
$\Spin(q_1) \times \ldots \times \Spin(q_d)$ acts on $C^0(q)_{\widetilde{E}}$ by
sending $(v_1, \ldots, v_d)$ to the endomorphism of $C^0(q)_{\widetilde{E}}$ given by
left multiplication with $m(v_1, \ldots, v_d) =  v_1 \cdot \ldots \cdot v_d$. Under this action each
$C^{\mathbf{a}}(q)$ is $(\Spin(q_1)\times \ldots \times \Spin(q_d))$-stable. Thus, by (i)
the decomposition (\ref{C0=PlusDa}) is also a decomposition of $\SMT(V)(\widetilde{E})$-modules.
Hence, by passing to $G$-invariants, (\ref{C0=PlusDa}) leads to a decomposition of Hodge
structures.

\vspace{2ex}
Denote by 
\begin{equation*}
  R:= D^{\mathbf{a}_0} = C^{\mathbf{a}_0}(q) = C^0(q_1) \otimes_{\widetilde{E}} \ldots \otimes_{\widetilde{E}} C^0(q_d).
\end{equation*}
By Lemma \ref{twistedC(Q)}, using the notations of Section \ref{SectionCorestriction},
we have 
\begin{equation*} 
 R = Z_G(C^0(Q))
\end{equation*}
as $G$-modules and hence $R^G = \cores_{E/\QQ} (C^0(Q))$. Thus we have recovered 
\begin{equation*}
W= \cores_{E/\QQ} (C^0(Q)) \subset C^0(q) = V
\end{equation*}
as a sub-Hodge structure. We now prove that after passing to $G$-invariants,
the remaining summands in (\ref{C0=PlusDa}) are
isomorphic to sums of copies of $W$.
\vspace{2ex}

Denote by $d_j = \sharp (G/ G_{\mathbf{a}_j})$ and 
choose a set of representatives $\mu_1, \ldots , \mu_{d_j}$ of $G/G_{\mathbf{a}_j}$ in $G$. 
We consider three group actions on $R^{\oplus d_j}$:

\vspace{1ex}
\hspace{2ex} $\bullet$ First there is a natural $(\Spin(q_1) \times \ldots \times \Spin(q_d))$-action which is just the diagonal
action of the one on $R$.

\vspace{1ex}
\hspace{2ex} $\bullet$ Let $\alpha : G \times R^{\oplus d_j} \to R^{\oplus d_j}$ be the diagonal action of the $G$-action
on $R$. 

\vspace{1ex}
\hspace{2ex} $\bullet$ Finally define the $G$-action $\beta$ by
 \begin{equation*} \beta: \left\{ \begin{aligned} 
 G \times \bigoplus_{l=1}^{d_j} R_{[\mu_l]} \to & \;
 \bigoplus_{l=1}^{d_j} R_{[\mu_l]} \\
 (\tau, (r_{[\mu_1]}, \ldots, r_{[\mu_d}])) \mapsto & \; (\tau r_{[\tau^{-1} \mu_1]}, \ldots, 
 \tau r_{[\tau^{-1} \mu_{d_j}]}). \end{aligned} \right.
\end{equation*} 

\vspace{2ex}
Now we will proceed in two steps:

\vspace{1ex}
(a) We show that $D^{\mathbf{a}_j}$ is isomorphic as $G$-module and as $(\Spin(q_1) \times
\ldots \times \Spin(q_d))$-module to $R^{\oplus d_j}$ where $G$ acts on the latter via $\beta$.

\vspace{1ex}
(b) We show that $R^{\oplus d_j}$ is isomorphic as $G$-module and as $(\Spin(q_1) \times 
\ldots \times \Spin(q_d))$-module with $G$ acting via $\alpha$ to $R^{\oplus d_j}$ with 
$G$ acting via $\beta$.

\vspace{2ex}
Note that neither of these two isomorphisms is canonical.
Once (a) and (b) are proved, we have an isomorphism
\begin{equation*}
 V_{\widetilde{E}} = C^0(q)_{\widetilde{E}} \simeq R^{\oplus 2^{d-1}}
\end{equation*}
 of $G$-modules and of $\SMT(V)(\widetilde{E})$-modules, $G$ acting diagonally on the right hand
 side. Here we use that
 \begin{equation*}
 \sum_j d_j = \sharp \left\{ \mathbf{a} \in \{ 0,1\}^d \; | \; \sum_i a_i \equiv 0 \; (2) \right\} = 2^{d-1}.
 \end{equation*}
 The proof of (ii) is then accomplished by passing to $G$-invariants.

\vspace{2ex} \noindent \it Proof of (a). \rm Denote by $F_j$ the field $\widetilde{E}^{G_{\mathbf{a_j}}}$.
As $C^{\mathbf{a}_j}(q) \subset D^{\mathbf{a}_j}$ is $G_{\mathbf{a}_j}$-stable,
$C^{\mathbf{a}_j}(q) = W_j \otimes_{F_j} \widetilde{E}$ for some $F_j$-vector
space $W_j$. Since $C^{\mathbf{a}_j}$ contains units in $C(q)_{\widetilde{E}}$, 
so does $W_j \subset C^{\mathbf{a}_j}$. (Very formally: There is a linear map $C^{\mathbf{a}_j} \to
\End(C(q)_{\widetilde{E}}), \; w \mapsto \{ v \mapsto v \cdot w \}$ which is defined over $F_{j}$.
The image of this map over $\widetilde{E}$ intersects the Zariski-open subset of automorphisms of
$C(q)_{\widetilde{E}}$, hence this must happen already over $F_j$.)

Choose a unit $w_j \in W_j$. Then for $\tau \in G$, since $w_j$ is $G_{\mathbf{a}_j}$-invariant, 
$\tau w_j \in C^{\tau {\mathbf{a}_j}}(q)$ depends only on the
coset $\tau G_{\mathbf{a}_j}$ and is again a unit in $C(q)_{\widetilde{E}}$.

Define an isomorphism of $\widetilde{E}$-vector spaces
\begin{equation*} \varphi: \left\{ \begin{aligned} 
 D^{\mathbf{a}_j} = \bigoplus_{l=1}^{d_j} C^{\mu_l \mathbf{a}_j} (q) \to
 & \; \bigoplus_{l=1}^{d_j} R_{[\mu_l]} \\
 (v_{\mu_1}, \ldots , v_{\mu_{d_j}}) \mapsto & \; (v_{\mu_1} \cdot \mu_1(w_j), \ldots , v_{\mu_{d_j}} \cdot \mu_{d_j} (w_j)).
 \end{aligned} \right.
\end{equation*}
This map is clearly $(\Spin(q_1) \times \ldots \times \Spin(q_d))$-equivariant since this group acts by multiplication
on the left whereas we multiply on the right.

As for the $G$-equivariance ($G$ acting via $\beta$ on the right hand side), we find
for \\ $(v_{[\mu_1]}, \ldots , v_{[\mu_{d_j}]}) \in D^{\mathbf{a}_j}$ and $\tau \in G$:
\begin{equation*} \begin{aligned}
 \varphi \big( \tau(v_{[\mu_1]}, \ldots,  v_{[\mu_{d_j}]}) \big) = & \;
 \varphi \big( \tau v_{[\tau^{-1} \mu_1]}, \ldots , \tau v_{[\tau^{-1} \mu_d]} \big) \\
 = & \; \big( \tau v_{[\tau^{-1} \mu_1]} \cdot \mu_1 w_j , \ldots , \tau v_{[\tau^{-1} \mu_{d_j}]} 
 \cdot \mu_{d_j} w_j \big) \\
 = & \; \big( \tau ( v_{[\tau^{-1} \mu_1]} \cdot \tau^{-1} \mu_1 w_j ), \ldots, 
      \tau ( v_{[\tau^{-1} \mu_{d_j}]} \cdot \tau^{-1} \mu_{d_j} w_j ) \big) \\
 = & \; \beta \big( \tau, ( v_{\mu_1} \cdot \mu_1 w_j, \ldots , v_{\mu_d} \cdot \mu_{d_j} w_j ) \big) \\
 = & \; \beta \big( \tau, \varphi(v_{[\mu_1]}, \ldots , v_{[\mu_{d_j}]}) \big).
 \end{aligned}
\end{equation*}
Here we used in the penultimate equality that $\sigma w_j$ depends only on the coset $\sigma G_{\mathbf{a}_j}$.
This proves (a). \hfill  $\mathrm{(a)}  \Box$ \par

\vspace{2ex}

\noindent \it Proof of (b). \rm Choose a $\QQ$-basis $f_1, \ldots, f_{d_j}$ of $F_j$. For $i=1, \ldots, d_j$ define
an $\widetilde{E}$-vector space homomorphism by
\begin{equation*}  \psi_i: \left\{ \begin{aligned}
  R \hookrightarrow & \; \bigoplus_{l=1}^{d_j} R_{[\mu_l]} \\
  r \mapsto & \; \big( \mu_1 (f_i) \cdot r , \ldots , \mu_{d_j} (f_i) \cdot r \big). \end{aligned} \right.
\end{equation*}
As $(\Spin(q_1) \times \ldots \times \Spin(q_d))(\widetilde{E})$ acts by $\widetilde{E}$-linear
automorphisms on $R$, the $\psi_i$ are equivariant for the Spin-action.

Let's show that $\psi_i$ is $G$-equivariant, $G$ acting on the right hand side via $\beta$. For $\tau \in G$ and
$r \in R$ we get
\begin{equation*} \begin{aligned}
 \psi_i( \tau r) = & \;\big(\mu_1 (f_i) \cdot \tau r, \ldots, \mu_{d_j} (f_i) \cdot \tau r \big) \\
    = & \; \big( \tau (\tau^{-1} \mu_1 (f_i) \cdot r) , \ldots , \tau ( \tau^{-1} \mu_{d_j}(f_i) \cdot r) \big) \\
    = & \; \beta \big( \tau, (\mu_1 (f_i) \cdot r, \ldots, \mu_{d_j} (f_i) \cdot r) \big) \\
    = & \; \beta(\tau, \psi_i (r)).
 \end{aligned}
\end{equation*}
Once more, we used the fact that $\sigma f_i$ depends only on the coset $\sigma G_{\mathbf{a}_j}$.

Finally, using Artin's independence of characters (see \cite[Thm.\ VI.4.1]{La}), we get
\begin{equation*}
  \det( (\mu_l (f_i))_{l,i}) \neq 0.
\end{equation*}
Consequently, the map
\begin{equation*}
 \oplus_{i=1}^{d_j} \psi_i : R^{\oplus d_j} \to R^{\oplus d_j}
\end{equation*}
is an isomorphism which has the equivariance properties we want and (b) is proved. \hfill $\text{(ii)} \Box$

\vspace{2ex}
\noindent
\it Proof of (iii). \rm Using that endomorphisms of Hodge structures are precisely those endomorphisms
which commute with the special Mumford--Tate group, we have to show that 
\begin{equation*}
 \End_{\SMT(V)} (W) = \cores_{E/\QQ}(C^0(Q)).
\end{equation*}
Denote
by $\mathfrak{g}$ the Lie algebra of $\SMT(V)$. Then 
\begin{equation*} \begin{aligned} 
 \End_{\SMT(V)} (W) = & \; \End_{\mathfrak{g}}(W) \\
  = & \; \{ f \in \End_{\QQ}(W) \; | \; Xf - fX = 0 \; \text{for all} \; X \in \mathfrak{g} \}.
  \end{aligned}
\end{equation*}
Since for any field extension $K /\QQ$ we have $\mathrm{Lie}(\SMT(V)_K) = \mathfrak{g} \otimes_\QQ K$
this implies that
\begin{equation} \label{TensK}
 \End_{\SMT(V)_K} (W_K) = \End_{\SMT(V)} (W) \otimes_{\QQ} K.
\end{equation}

Now $\SMT(V)(\widetilde{E}) = \Spin(q_1) \times \ldots \times \Spin(q_d) (\widetilde{E})$ acts on 
$W_{\widetilde{E}} = C^0(q_1) \otimes \ldots \otimes C^0(q_d)$ by factorwise left multiplication:
\begin{equation*}
 \big( (v_1, \ldots, v_d), w_1 \otimes \ldots \otimes w_d \big) \mapsto
  (v_1 \cdot w_1) \otimes \ldots \otimes (v_d \cdot w_d).
\end{equation*}
Therefore, using multiplication on the right, we get an inclusion
\begin{equation*}
  \big( C^0(q_1) \otimes \ldots \otimes C^0(q_d) \big)^{\mathrm{op}} \hookrightarrow 
  	\End_{\SMT(V)(\widetilde{E})} (W_{\widetilde{E}}), \; \; \;
  w \mapsto \{ w' \mapsto w' \cdot w \}. 
\end{equation*}
Now, $( C^0(q_1) \otimes \ldots \otimes C^0(q_d) )^{\mathrm{op}} \simeq C^0(q_1)^{\mathrm{op}} 
\otimes \ldots \otimes C^0(q_d)^{\mathrm{op}} \simeq C^0(q_1) \otimes \ldots \otimes C^0(q_d)$
and hence passing to $G$-invariants we have an inclusion
\begin{equation} \label{inclusion}
 \cores_{E /\QQ} (C^0(Q)) \hookrightarrow \End_{\SMT(V)(\QQ)} (W).
\end{equation}

We will now show that this is an isomorphism over $\widetilde{E}$. Using (\ref{TensK}) and comparing dimensions this
will prove (iii).

To show that (\ref{inclusion}) is an isomorphism over $\widetilde{E}$ we have to 
determine the $\Spin(q_1) \times \ldots \times \Spin(q_d)$-invariants
in 
\begin{equation*}
  \End_{\widetilde{E}}\left( C^0(q_{1}) \otimes \ldots \otimes C^0(q_d) \right)
= \End_{\widetilde{E}}C^0(q_1) \otimes \ldots \otimes \End_{\widetilde{E}}C^0(q_d).
\end{equation*}
Using the next lemma inductively, this is equal to
\begin{equation*}
 \End_{\Spin(q_1)} C^0(q_{1}) \otimes \ldots \otimes \End_{\Spin(q_d)} C^0(q_{d}).
\end{equation*}
Now by \cite[Lemma 6.5]{vG1}, $\End_{\Spin(q_i)} C^0(q_{i}) = C^0(q_{i})$. This proves (iii). \hfill $\Box$

\begin{Lemma} Let $G$ and $H$ be two reductive linear algebraic groups over a field $K$ of charac\-teristic $0$. 
 Let $M$ resp.\ $N$  be finite-dimensional representations over $K$ of $G$ resp.\ $H$. Then 
\begin{equation*}
  (M \otimes_{K} N)^{G \times H} = M^G \otimes_{K} N^H.
\end{equation*}
\end{Lemma}

\begin{Proof}
 Decompose $M = \bigoplus_i M_i$ and $N = \bigoplus_j N_j$ in irreducible representations. Then
 $M_i \otimes N_j$ is an irreducible representation of $G \times H$ since fixing $0 \neq m_0 \in M_i$
 and $0 \neq n_0 \in N_i$ the orbit $(G \times H) m_0 \otimes n_0$ generates $M_i \otimes N_j$.
 
 To conclude the proof note that the space of invariants is the direct sum of trivial, one-dimensional
 sub representations.
\end{Proof}
 
\end{subsection}

\begin{subsection}{The Brauer--Hasse--Noether theorem}
Let $k$ be a field of characteristic $\neq 2$, let $A$ be a central simple $k$-algebra (i.e.\ a finite-dimensional
$k$-algebra with center $k$ which has no non-trivial two-sided ideals). By Wedderburn's theorem, there
exists a central division algebra $D$ over $k$ and an integer $n > 0$ such that $A \simeq \mathrm{Mat}_n(D)$.
Let $d^2$ be the dimension of $D$ over $k$ (this is a square because after base change, $D$ becomes isomorphic
to a matrix algebra). Then $d$ is the \emph{index of $A$}, denoted by $i(A)$. 
The class of $A$ in the Brauer group of $k$ has finite order. This integer is called the
\emph{exponent of $A$}, it is denoted by $e(A)$. In general, we have $e(A) | i(A)$.

\vspace{2ex}
Let $K / k$ be a cyclic extension of degree $n$, let $\sigma$ be a generator of the Galois group
$\mathrm{Gal}(K/k)$, let $a \in k^*$. There is a central simple $k$-algebra $(\sigma, a, K/k)$ which
as a $k$-algebra is generated by $K$ and an element $y \in (\sigma, a, K/k)$ such that
\begin{equation*}
 y^n = a \; \; \text{and} \; \; r \cdot  y =  y \cdot \sigma(r) \; \text{for} \; r \in K.
\end{equation*}
This algebra is called the \emph{cyclic algebra associated with $\sigma, a$ and $K/k$}.
A cyclic algebra over $k$ of dimension 4 is a quaternion algebra.

\begin{Thm}[Brauer, Hasse, Noether \cite{BHN}] \label{BrauerHasseNoether}
 Let $k$ be an algebraic number field. Then any central division algebra $A$ over $k$ is a cyclic algebra (for an appropriate 
 cyclic extension $K/k$ and $\sigma$ and $a$ as above). Moreover, the exponent and the index of $A$ coincide. In particular,
 a central division algebra of exponent 2 is a quaternion algebra.
\end{Thm}

\end{subsection}

\begin{subsection}{An example}
 We continue to assume that $(T,h,q)$ is a Hodge structure 
 of K3 type with $E = \End_{\Hdg}(T)$ a totally real number field 
 of degree $d$ over $\QQ$. By \cite[Prop.\ 3.2]{vG2} we have $\dim_E T \ge 3$. We will consider now the case that
 $\dim_E T = 3$.

 Then $T_1$ is an $3$-dimensional $\widetilde{E}$-vector space with quadratic form $q_1$ of
 signature $(2+, 1-)$. The $3$-dimensional quadratic spaces $(T_2,q_2), \ldots, (T_d,q_d)$ are negative definite.
 This implies that
 \begin{equation*} \begin{aligned}
   C^0(q_1)_{\RR} & = \mathrm{Mat}_2 (\RR) \; \text{and} \\
   C^0(q_i)_{\RR}&  = \HH \; \text{for} \; i \ge 2 \end{aligned}
 \end{equation*}
 (see \cite[Thm.\ 7.7]{vG1}).
 Since 
 \begin{equation*}
  \cores_{E/\QQ} (C^0(Q)) \otimes_{\QQ} \widetilde{E} = Z_G(C^0(Q)) = C^0(q_1) \otimes_{\widetilde{E}} \ldots \otimes_{\widetilde{E}} C^0(q_d)
 \end{equation*} 
 we get
 \begin{equation*}
  \cores_{E/\QQ} (C^0(Q)) \otimes_{\QQ} \RR = \mathrm{Mat}_{2} (\RR) \otimes_{\RR} \HH \otimes_{\RR} \ldots \otimes_{\RR} 
   \HH.
 \end{equation*}
 Now, since $\HH \otimes \HH \simeq \mathrm{Mat}_4(\RR)$ this becomes 
 \begin{equation} \label{CoresTensR}
   \cores_{E/\QQ} (C^0(Q)) \otimes_{\QQ} \RR \simeq \left\{ 
    \begin{aligned} & \mathrm{Mat}_{2^{d-1}} (\HH) \; \text{for even} \; d \\
      & \mathrm{Mat}_{2^{d}} (\RR) \; \text{for odd} \; d. 
    \end{aligned} \right.
 \end{equation}
 On the other hand, the corestriction induces a homomorphism of Brauer groups
 \begin{equation*}
  \mathrm{cores} : \mathrm{Br}(E) \to \mathrm{Br}(\QQ)
 \end{equation*}
 (cf.\ \cite[\S 9, Thm.\ 5]{D}). Therefore, the exponent of $\cores_{E/\QQ}(C^0(Q))$ in the Brauer
 group of $\QQ$ is 2.
 By the Brauer--Hasse--Noether Theorem \ref{BrauerHasseNoether} there exists a (possibly split) quaternion algebra $D$ over $\QQ$ with
 \begin{equation} \label{Cores=D}
   \cores_{E/\QQ} (C^0(Q)) \simeq \mathrm{Mat}_{2^{d-1}} (D).
 \end{equation}
 Combining (\ref{CoresTensR}) with (\ref{Cores=D}) we see that $D$ is a definite quaternion algebra over $\QQ$
 in case $d$ is even and an indefinite quaternion algebra in case $d$ is odd.
 The endomorphism algebra of a Kuga--Satake variety of $(T,h,q)$ is $\mathrm{Mat}_{2^{2d-2}}(D)$. 
 Since the dimension of a Kuga--Satake variety is $2^{\dim_{\QQ} (T) -2} = 2^{3d-2}$, we have proved
 \begin{Cor} \label{KorRM}
   Let $(T,q,h)$ be a Hodge structure of K3 type with $E = \End_{\Hdg}(T)$ a totally real number field of degree $d$ over $\QQ$. Assume that
   $\dim_E(T)=3$. Then for any Kuga--Satake variety $A$ of $(T,h,q)$ there exists an isogeny
   \begin{equation*}
     A \sim B^{2^{2d-2}}
   \end{equation*}
   where $B$ is a $2^d$-dimensional Abelian variety.

   \hspace{1ex} If $d$ is even, $B$ is a simple Abelian variety of type III, i.e. $\End_{\QQ} (B) = D$ for a definite quaternion algebra $D$
   over $\QQ$.

   \hspace{1ex} If $d$ is odd, $B$ has endomorphism algebra $\End_{\QQ}(B) =D$ for an indefinite (possibly split) quaternion algebra $D$
   over $\QQ$.
 \end{Cor} 
 
 \noindent
 {\it Remark.} (i) In the case $d = 2$ and $\dim_E(T) =3$, van Geemen showed in \cite[Prop.\ 5.7]{vG2} that the Kuga--Satake 
  variety of $T$ is isogenous to a self-product of an Abelian fourfold with definite quaternion multiplication
  and Picard number 1. It is this case which
  will be of interest in the next section.

  \vspace{1ex}
  (ii) The case $d=\dim_E (T) =3$ was also treated by van Geemen (see \cite[5.8 and 6.4]{vG2}). 
  He considers the case $D \simeq \mathrm{Mat}_2(\QQ)$ and
  relates this to work of Mumford and Galluzzi. Note that in this case the Abelian variety $B$ of the corollary
  is not simple.

  \vspace{2ex}
  \noindent
  {\it Example.} In \cite[3.4]{vG2}, van Geemen constructs a one-dimensional 
  family of six-dimensional K3 type Hodge structures with
  real multiplication by a quadratic field $E= \QQ(\sqrt{d})$ for some square-free integer $d> 0$ which 
  can be written in the form $d = c^2 + e^2$ for rational $c,e>0$. 
  These Hodge structures 
  are realized as the transcendental lattice of certain K3 surfaces which are double 
  covers of $\PP^2$, see Section \ref{SectionDoubleCovers}.
  Pick a member $S$ of this family. Then $T(S) \otimes_{\QQ} E$ splits in the direct sum of 
  two three-dimensional $E$-vector spaces $T_1$ and 
  $T_2$. It turns out that the quadratic space $(T_1, q_1) = (T_1,Q)$ is isometric to 
  $(E^3, \sqrt{d} X_1^2 + \sqrt{d} X_2^2 - (d-\sqrt{d}c) X_3^2)$. Consequently
  \begin{equation*}
    C^0(Q) = (-d, \sqrt{d} (d - \sqrt{d}c))_E \simeq (-1, \sqrt{d}- c)_E.
  \end{equation*}
  Here for $a,b \in E^*$, the symbol $(a,b)_E$ denotes the quaternion algebra over $E$ generated by elements 
  $1, i$ and $j$ subject to the relations $i^2 = a, j^2 =b$ and $ij = -ji$ (see \cite[Ex.\ 7.5]{vG1}).

  The projection formula for central simple algebras (see \cite[Thm.\ 3.2]{T}) implies that
  \begin{equation*} \begin{aligned}
   \cores_{E/\QQ} (C^0(\QQ)) & \simeq (-1, N_{E/\QQ}( \sqrt{d} - c))_{\QQ} \\
     & \simeq (-1, c^2 - d)_{\QQ} \simeq (-1, - e^2)_{\QQ} \simeq (-1, -1)_{\QQ} \end{aligned}
  \end{equation*}
  which are simply Hamilton's quaternions over $\QQ$. Here, $N_{E/\QQ}: E \to \QQ$ is
  the norm map.
  Hence, a Kuga--Satake variety for $T(S)$ is isogenous to a 
  self-product $B^4$ where $B$ is a simple Abelian fourfold
  with $\End_{\QQ}(B) = (-1, -1)_{\QQ}$.

\end{subsection}
\end{section}

\begin{section}{Double covers of \texorpdfstring{$\PP^2$}{P2} branched along six lines} \label{SectionDoubleCovers}
Let $S$ be a K3 surface which admits a morphism $p: S \to \PP^2$ such that the branch locus of $p$ is 
the union of six lines. 

In this section we use the decomposition theorem to prove Theorem \ref{ThmHCDC} which states that
the Hodge conjecture holds for $S \times S$.

 \begin{subsection}{Abelian varieties of Weil type} \label{AVWeilType}
 By a result of Lombardo \cite{L}, the Kuga--Satake variety of $S$ is
 of Weil type. We briefly recall what this means.
 
 \vspace{2ex}
 Let $K = \QQ(\sqrt{-d})$ for some square-free $d \in \NN$. A polarized Abelian variety $(A,H)$
 of dimension $2n$ 
 is said to be of \emph{Weil type for $K$} if there is an inclusion $K \subset \End_{\QQ}(A)$ 
 mapping $\sqrt{-d}$ to $\varphi$ such that
 
 \vspace{1ex}
 \hspace{2ex} $\bullet$ the restriction of $\varphi^*: H^1(A,\CC) \to H^1(A,\CC)$ to $H^{1,0}(A)$ is diagonalizable with eigenvalues     
 $\sqrt{-d}$ and $- \sqrt{-d}$, both appearing with multiplicity $n$,
 
 \vspace{1ex}
 \hspace{2ex} $\bullet$ $\varphi^* H = d H$.
 
 \vspace{1ex}
 
 There is a natural $K$-valued Hermitian form on the $K$-vector space $H^1(A,\QQ)$ which is defined
 by 
 \begin{equation*} \widetilde{H}:   H^1(A,\QQ) \times H^1(A,\QQ) \to K, \; \; \;
  (v,w) \mapsto H(\varphi^* v, w) + \sqrt{-d} H (v,w).
 \end{equation*}
 By definition, the discriminant of a polarized Abelian variety of Weil type $(A,H,K)$ is
 \begin{equation*}
  \mathrm{disc}(A,H,K) = \mathrm{disc} (\widetilde{H}) \in \QQ^* / N_{K/\QQ}(K^*)
 \end{equation*}
 where $N_{K/\QQ}: K \to \QQ$ is the norm map.
 
 Polarized Abelian varieties of Weil type come in $n^2$-dimensional families (see \cite[5.3]{vG3}).
  
 \vspace{2ex}
 Weil introduced such varieties as examples of Abelian varieties
 which carry interesting Hodge classes. He constructs a two-dimensional space, called 
 the space of \emph{Weil cycles}
 \begin{equation*}
  W_K \subset H^{n,n}(A,\QQ).
 \end{equation*}
 For the definition of $W_K$ see
 \cite[5.2]{vG3}. In general, the algebraicity of the classes 
 in $W_K$ is not known. 
 Nonetheless there are some positive results. Here we mention one which we will use below.
 \begin{Thm}[Schoen \cite{S} and van Geemen \cite{vG4}, Thm.\ 3.7] \label{ThmvG}
 Let $(A,H)$ be a polarized
 Abelian fourfold of Weil type for the
 field $\QQ(i)$. Assume that the discriminant of $(A,H,\QQ(i))$ is $1$. Then the space of Weil cycles
 $W_{\QQ(i)}$ is spanned by classes of algebraic cycles.
 \end{Thm}
  
 Van Geemen uses a six-dimensional eigenspace in the complete linear system of the unique totally symmetric line bundle
 $\mathcal{L}$ with $\cc_1(\mathcal{L}) = H$ 
to get a rational (2:1) map of $A$ onto a quadric $Q \subset \PP^5$. Then the projection on $W_{\QQ(i)}$ of the
 classes of the pullbacks of the two rulings
 of $Q$ generate the space $W_{\QQ(i)}$.
 \end{subsection}
 
 \begin{subsection}{Abelian varieties with quaternion multiplication}
 Let $D$ be a definite quaternion algebra over $\QQ$. Such a $D$ admits an involution $x \mapsto \overline{x}$
 which after tensoring with $\RR$ becomes the natural involution on Hamilton's quaternions $\HH$.

 \vspace{2ex} 
 A polarized Abelian variety $(A,H)$ of dimension $2n$ has \emph{quaternion multiplication by $D$}
 if there is an inclusion
 $D \subset \End_{\QQ}(A)$ such that 
 
 \vspace{1ex} 
 \hspace{2ex} $\bullet$ $H^1(A,\QQ)$ becomes a $D$-vector space and
 
 \vspace{1ex} 
 \hspace{2ex} $\bullet$ for $x \in D$ we have $x^* H = x \overline{x} H$.
 
 \vspace{1ex}
 We say that $(A,H,D)$ is an Abelian variety of definite quaternion type.
 Polarized Abelian varieties of dimension $2n$ 
 with quaternion multiplication by the same quaternion algebra come
 in $n (n-1)/2$-dimensional families (cf.\ \cite[Sect.\ 9.5]{BL}). 
 
 \vspace{2ex}
 Let $K \subset D$ be a quadratic extension field of $\QQ$. Then $K$ is a CM field and $(A,H,K)$ is a
 polarized Abelian variety of Weil type (see \cite[Lemma 4.5]{vGV}). The space of quaternion
 Weil cycles of $(A,H,D)$
 \begin{equation*}
   W_D \subset H^{n,n}(A,\QQ)
 \end{equation*}
 is defined as the span of $x^* W_K$ where $x$ runs over $D$. It can be shown that this 
 is independent of the choice of $K$ (see \cite[Prop.\ 4.7]{vGV}). For the general member of the
 family of polarized Abelian varieties 
 with quaternion multiplication these are essentially all Hodge classes:
 
 \begin{Thm}[Abdulali, see \cite{Ab}, Thm.\ 4.1] \label{ThmA}
 Let $(A,H,D)$ be a general Abelian variety 
 of quaternion type. Then the space of Hodge classes on any self-product of $A$ is generated by 
 products of divisor classes and quaternion Weil cycles on $A$. 
 
 In particular, if for one quadratic extension field $K \subset D$ the space of Weil cycles
 $W_K$ is known to be algebraic, then the Hodge conjecture holds for any self-product of $A$.
 \end{Thm}
 
 In Abdulali's theorem, a triple $(A,H,D)$ is general if the special Mumford--Tate group
 of $H^1(A,\QQ)$ is the maximal one. In the moduli space of triples $(A,H,D)$ the locus
 of general triples
 is everything but a countable union of proper, closed subsets.
 \end{subsection}
 
 \begin{subsection}{The transcendental lattice of \texorpdfstring{$S$}{S}} We now turn back to our
 K3 surface $S$.
Let $p: S \to \PP^2$ be the (2:1) morphism which is ramified over six lines.

The N\'eron--Severi group of $S$ contains the 15 classes $e_1, \ldots, e_{15}$ corres\-ponding to the 
exceptional divisors over the intersection points of the six lines. Let $h$ be the class of the pullback
of $\OO_{\PP^2}(1)$. 

Define $\widetilde{T}(S):= \langle e_1, \ldots, e_{15} ,h \rangle^{\perp} \subset H^2(S ,\QQ)$.
The (rational) transcendental lattice of $S$ is defined to be $T(S): = \NS(S)^{\perp} \subset
H^2(S,\QQ)$. Then we have
\begin{equation*}
  T(S) \subset \widetilde{T}(S).
\end{equation*}
Both, $T(S)$ 
and $\widetilde{T}(S)$ are Hodge structures of K3 type. In addition,
$T(S)$ is irreducible. Since the second Betti number of $S$ is 22, the $\QQ$-dimension of $\widetilde{T}(S)$
is 6.
\end{subsection}

 \begin{subsection}{The Kuga--Satake variety of \texorpdfstring{$\widetilde{T}(S)$}{T(S)}}
 Denote by $A$ the Kuga--Satake variety associated with $\widetilde{T}(S)$.
 \begin{Thm}[Lombardo, see \cite{L}, Cor.\ 6.3 and Thm.\ 6.4] \label{ThmLombardo} There is an isogeny
 \begin{equation*}
    A \sim B^4
 \end{equation*}
 where $B$ is an Abelian fourfold with $\QQ(i) \subset \End_{\QQ}(B)$. Moreover, 
 $B$ admits a polarization $H$ such that $(B,H,\QQ(i))$ is a polarized Abelian variety of Weil type
 with $\mathrm{disc}(B,H,\QQ(i))=1$. 
 \end{Thm}

 Paranjape \cite{P}
 explains in a very nice way how this variety $B$ is geometrically related to $S$.
 He shows that there exists a triple 
\begin{equation*} 
 (C, E, f: C \to E)
\end{equation*}
where $C$ is a genus five curve, $E$ an elliptic curve and $f$ a $(4:1)$ map 
such that 
\begin{equation*}
  \mathrm{Prym}(f) = B.
\end{equation*}
Then $S$ can be obtained as the resolution
of a certain quotient of $C \times C$. It is noteworthy that Paranjape does not
construct explicitly a triple $(C,E,f)$ starting with a K3 surface $S$ in 
the family $\pi$. His proof goes the other way round. He associates 
to any triple a K3 surface and shows then that letting vary the triple he obtains all surfaces 
in the family $\pi$.

Paranjape's construction establishes that the Kuga--Satake inclusion
\begin{equation} \label{Paranjape}
 \widetilde{T}(S) \hookrightarrow H^2(B^4 \times B^4, \QQ)
\end{equation}
is given by an algebraic cycle on $S \times B^4 \times B^4$.
\end{subsection}

\begin{subsection}{Proof of Theorem 2} As pointed out in the introduction, we
 have to prove that $E_S := \End_{\Hdg}(T(S))$ is spanned by algebraic classes.
 Since the Picard number of $S$ is at least 16, 
 we can apply Ram\'on-Mar\'i's corollary \cite{Ma} of Mukai's theorem \cite{Mu} which
 proves the assertion in the case that $S$ has complex multiplication.
 
 Therefore, we may assume that 
 $S$ has real multiplication. Note that $T(S)$ is an $E_S$-vector space and that
 $\dim_{E_S} T(S) \cdot [E_S :\QQ] = \dim_{\QQ} T(S) \le 6$. On the other hand, by \cite[Lemma 3.2]{vG2}, we know
 that $\dim_{E_S} T(S) \ge 3$. It follows that
 either $E_S =\QQ$ or $E_S = \QQ(\sqrt{d})$ for some square-free $d \in \QQ_{>0}$. In the first case
 we use the fact, that the class of the diagonal $\Delta \subset S \times S$ induces the identity on the cohomology
 and that the K\"unneth projectors are algebraic on surfaces
 so that $\QQ \id \subset E_S$ is spanned by an algebraic class. 

 It remains to study the case $E_S = \QQ(\sqrt{d})$. The idea is to consider
 the Kuga--Satake variety $A(S)$ of $\widetilde{T}(S) = T(S)$. By Paranjape's theorem the inclusion
 \begin{equation*}
   \widetilde{T}(S) \subset H^2(A(S) \times A(S),\QQ)
 \end{equation*}
 is algebraic. It follows that there is an algebraic projection
 $\pi: H^2(A(S) \times A(S),\QQ) \to \widetilde{T}(S)$ (see \cite[Cor.\ 3.14]{Kl}) and therefore it
 is enough to show that there
 is an algebraic class 
 \begin{equation*} 
   \alpha \in H^2(A(S) \times A(S),\QQ) \otimes H^2(A(S) \times A(S),\QQ)  \subset H^4(A(S)^4, \QQ)
 \end{equation*}
 with $\pi \otimes \pi (\alpha) = \sqrt{d}$.

 Combining
 Corollary \ref{KorRM} with Lombardo's theorem \ref{ThmLombardo} we see that
 $A(S) \sim B^4$ where $B$ is an Abelian fourfold with 
 $\End_{\QQ}(B) = D$ for a definite quaternion algebra
 and $\QQ(i) \subset D$. 
 Moreover, there is a polarization $H$
 of $B$ such that $(B,H,\QQ(i))$ is a polarized Abelian variety of Weil type
 of discriminant 1. Since by \cite[Prop.\ 5.5.7]{BL}, the
 Picard number of $B$ is 1, $(B,H,D)$ is a polarized
 Abelian variety of quaternion type.
 
 There is a one-dimensional family $(B,H,D)_t$ of deformations of $(B,H,D)$ and this corres\-ponds
 to a one-dimensional family $S_t$ of deformations of $S$ which parametrizes K3 surfaces with
 real multiplication by the same class. By Abdulali's Theorem \ref{ThmA}, for $t$ general
 the space of Hodge classes on $(B_t)^{16} \sim A(S_t)^4$ is generated by products of divisors and 
 quaternion Weil cycles, that is by products of $H$ and classes in $W_D$. Denote the span of these products in 
 $H^4(A(S_t)^{4},\QQ)$ by $F_t$.

 Since the class corresponding to $\sqrt{d} \in \widetilde{T}(S_t) \otimes \widetilde{T}(S_t)$,
 the projection $\pi: H^2(A(S_t)^2,\QQ) \to  \widetilde{T}(S_t)$ and the space $F_t$ are locally constant, 
 there $\mbox{exists}$ a locally constant class $\alpha_t \in H^4(A_{S_t},\QQ)$ with the properties:
 
 $\bullet$ for all $t$ we have $\pi \otimes \pi (\alpha_t) = \sqrt{d}$,
 
 $\bullet$ for all $t$ we have $\alpha_t  \in F_t$.
 
 Now by Schoen's and van Geemen's Theorem \ref{ThmvG} the space of Weil cycles $W_{\QQ(i)}$
 is generated by algebraic classes on any $B_t$. It follows that $W_D$ is generated by algebraic classes and consequently $F_t$
 is generated by algebraic classes for any $t$. In particular, $\alpha_t \in F_t$
 is algebraic. This proves the theorem. \hfill $\Box$ \par
\end{subsection}
\end{section}

\end{document}